\documentclass{article}
\usepackage{amsfonts,amsmath,amssymb,textcomp,enumerate,bbm}
\usepackage{amsthm}
\newcommand{\R}{\mathbb{R}}
\DeclareMathOperator{\Ric} {Ric}

\newtheorem{theorem}{Theorem}

\newtheorem{lemma}{Lemma}

\usepackage[english]{babel}

\newcommand{\assign}{:=}

\newcommand{\tmop}[1]{\ensuremath{\operatorname{#1}}}

\title{Dimensionless $L^p$ estimates for the Riesz vector on manifolds} 
\author{K. Dahmani, K. Domelevo, S. Petermichl} 
\date{}

\begin{document} 

\maketitle
\begin{abstract}
We present a new proof of the dimensionless $L^p$ boundedness of the Riesz vector on manifolds with bounded geometry. Our proof has the significant advantage that it allows for a much stronger conclusion, namely that of a new dimensionless weighted $L^p$ estimate with optimal exponent. 
Other than previous arguments, only a small part of our proof is based on special auxiliary functions, the core of the argument is a weak type estimate and a sparse decomposition of the stochastic process by X.D. Li, whose projection is the Riesz vector.
\end{abstract}

\section{Introduction}

In this paper, we are interested in dimensionless weighted and unweighted $L^p$ norm estimates of the Riesz vector on manifold. In the Euclidian setting, the $i$--th Riesz transfrom in $\R^n$ is defined as
\[
        R_i = \frac{\partial}{\partial x_i} (-\Delta)^{-1/2},
\]
where $\Delta = \sum_{i=1}^n \partial^2 / \partial_{x_i}^2$ is the usual Laplacian in $\R^n$. The vector Riesz transform $R$ is defined as the collection $R=(R_1,R_2,\ldots,R_n)$. In the one-dimensional setting, the Riesz transform is nothing but the Hilbert transform. The $L^p$ estimate of the Hilbert transform on the real line dates back to the work of Riesz \cite{Riesz1927} and Pichorides \cite{Pic1972}. Regarding the $L^p$ estimate of the Riesz vector in $\R^n$, see \cite{Stein1955, Meyer1984, Pisier1988, BanWan1995, IwaMar1996, DraVol2006}.

A corner stone in this line of results is the stochastic representation of Riesz transfroms in $\R^n$ by Gundy--Varopoulos \cite{GV1979}. To this end, these authors define the so-called background noise which are Brownian trajectories in the upper half space started at infinity and stopped when hitting the boundary. To a given function $f$ defined on $\R^n$ and its Poisson extension in the upper half space, these authors associate a natural martingale $M^f$. They prove that the Riesz transforms can be written as a suitable conditional expectation of martingale transforms of $M^f$. This representation was extended to the Riemannian manifold setting by X.-D. Li \cite{Li2008,LiErratum,LiArXiv} and has thus enabled the first dimensionless estimates with the growth proportional to $(p-1)^{-1}$ when $p>2$ and $p-1$ when $p>2$ in this setting \cite{CD}, \cite{BO}.

For early considerations of $L^p$ boundedness of Riesz transforms on manifolds, we refer to \cite{SteinTopics}. We mention also the works \cite{Meyer1984,Bakry1985,Gundy1986,Bakry1987,Pisier1988,Arcozzi} among which the papers of Bakry provide estimates of Riesz transforms for complete Riemannian manifolds under the general condition that the Bakry--Emery curvature is bounded below (see \cite{Emery}). Using stochastic techniques, linear dimensionless estimates of the Bakry--Riesz vector on manifold were announced in \cite{Li2008} \cite{LiErratum}. Using deterministic techniques, such estimates were proved in \cite{CD}. See also \cite{BanBau2013} for second order Riesz transforms on manifolds and \cite{BO} for Riesz transforms on manifolds, correcting a previous gap in the probabilistic proof.

In this paper, we will only consider manifolds with non-negative curvature. We will use stochastic tools relying on the stochastic representation of Riesz transforms on manifolds by X.-D. Li \cite{Li2008,LiErratum,LiArXiv}.

Our proof is very different from previous ones in that it does not rely on a Bellman function for the problem. Rather, it develops a sparse domination of the stochastic process of Li. See the elegant and short argument in \cite{Lacey2015} for the first probabilistic object, a discrete time martingale transform and also \cite{DP2} for the continuous time case. One can deduce, from such domination a dimensionless bound. The sparse operators are particularly well suited for working with weights, which is why this so obtained dimensionless estimate also holds in the weighted setting.

The stochastic process by Li is a specific semi-martingale, built using a pair of martingales that have differential subordination and solving a certain stochastic differential equation.  As such, our argument required several new tools. One of them is a weak type estimate of the maximal operator of this process. This is the only part of our proof that uses a (simple) Bellman function. The explicit form of the function is essential and not just its convexity and size properties. The first derivative of said Bellman function is used to control a drift term that arises because the process we consider is not a martingale. 
Further, we then show that this process has a sparse domination, according to the definition of sparse operator in  \cite{DP2}. The specific form of the defining stochastic equation is used.

The rest of the arguments, to deduce the estimates for the Riesz vector, are considered standard.

\subsection{Bakry--Riesz transforms on manifolds.}

Let $(M, g)$ be a complete Riemannian manifold with metric $g$ and dimension
$n$. Let $\Delta$ be the non-positive Laplace--Beltrami operator, $\nabla$
the gradient operator, and $\nu$ the volume measure on $(M, g)$ such that
$d  \nu (x) = \sqrt{\det g (x)} d  x$. Let moreover $\mu_{\varphi}$ be a
weighted volume measure on $M$ with $d  \mu_{\varphi} (x) = e^{- \varphi (x)} d 
\nu (x)$, where $\varphi (x) \in \mathcal{C}^2 (M)$. The weighted Laplacian
$\Delta_{\varphi}$ with respect to $\mu_{\varphi}$ on $M$ is defined for any function
$f$ by
\[ \Delta_{\varphi} f \assign \Delta - \nabla \varphi \cdot \nabla f. \]

We assume that $\mu_{\varphi} (M) < \infty$ and by normalizing, we may assume without
loss of generality that $\mu_{\varphi}$ is actually a probability measure. The
Bakry--Emery curvature tensor associated with $\Delta_{\varphi}$ is defined by
\[ \tmop{Ric}_{\varphi} = \tmop{Ric} + \nabla^2 \varphi, \]
where $\tmop{Ric}$ denotes the Ricci curvature tensor on $M$ and $\nabla^2
\varphi$ the Hessian of $\varphi$ with respect to the Levi-Civita connection
on $(M, g)$.

All over this paper, we assume a non-negative curvature
\[
	\tmop{Ric}_{\varphi}  \geq 0.
\]
We denote by $R$ the Bakry-Riesz transform defined as
\[
	R = d \circ (-\Delta _{\varphi} )^{-1/2}.
\]
It was proved by Bakry \cite{Bakry1987} that for any $p>1$ there exists a universal constant $C_p$ 
such that for any function $f\in C_0^\infty(M)$, there holds
\[
	\|R  f \|_{L^p(\mu)} \leq C_p \| f\|_{L^p(\mu)}.
\]

\paragraph{Probabilistic setting and notations.}
Let $B^M_t$ the diffusion process on $M$ with generator $\Delta_{\varphi}$ obeying
\[
	dB^M_t = U_t dW_t - \nabla\varphi(B^M_t) dt,
\]
where $W_t$ is the Brownian motion on $\R^n$, $U_t \in \tmop{End}(T_x M,T_x M)$ denotes the stochastic parallel transport on $M$ along the trajectory $\{  B^M_s  , 0\leq s \leq t\}$. Let further $B_t$ the one-dimensional Brownian motion started at $y>0$ with the normalisation $E [(B^M_t)^2] = 2t$ such that its generator is $\partial^2/\partial y^2$.
Following \cite{Meyer1984,GV1979, Gundy1986}, there exists a diffusion process $Z_t=(B^M_t,B_t)$ on $(M,\R^+)$ -- the so-called background radiation process -- associated to the generator $\Delta_\varphi + \partial^2/\partial y^2$ and with initial distribution $\mu \otimes \delta_y$.

We recall that the martingale $Y$ is said differentially subordinate to the martingale $X$ if the process $(\langle X,X\rangle_t-\langle Y,Y\rangle _t)_{t\geq0}$ is non-negative and non-decreasing in $t$, where the bracket $\langle \cdot,\cdot \rangle$ denotes the usual quadratic covariance process for real or vector-valued processes.

\paragraph{Probabilistic representation of the Bakry-Riesz vector on manifolds.}

Using a martingale approach, one can represent the Riesz vector $R$ via a probabilistic representation. In the literature, it first appeared in \cite{GV1979}, where the Riesz transform was defined on $\mathbb{R}^n$. In \cite{Arcozzi} Arcozzi extended this formula to compact Lie groups and spheres. In \cite{Li2008}, \cite{LiErratum} Li presented a new formula adapted to complete Riemannian manifolds. The representation formula of the Riesz vector in this setting for a complete manifold with $\Ric _{\varphi}\geq 0$ is as follows

\begin{equation} \label{rep}
-\frac{1}{2}(R f)(x)=\lim_{y\rightarrow \infty} \mathbb{E}_y\left[M_{\tau}\int_0^{\tau}  M_s^{-1} dQ(f)(B^M_s,B_s)dB_s | B^M_{\tau}=x \right],
\end{equation}
where 
\begin{itemize}
\item $Q(f)(x,y)=e^{-y\sqrt{-\Delta_{\varphi}}}f(x)$ is the Poisson semigroup;
\item $\tau = \inf \{t>0:B_t=0\}$ is the stopping time upon hitting the boundary of the upper half space;
\item $M_t$ is the solution to the matrix-valued stochastic differential equation
\[dM_t=V_tM_tdt, \ \ \ M_0=Id,\]
for some adapted and continuous process $(V_t)_{t \geq 0}$ taking values in the set of symmetric and non-positive $d \times d$ matrices.
\end{itemize} 

Equivalently, one can rewrite this fomula as
\begin{equation}
-\frac{1}{2}(R f)(x)=\lim_{y\rightarrow \infty} \mathbb{E}_y\left[   Z_\tau   |B^M_{\tau}=x \right],
\end{equation}
where $Z_t$ is a semi-martingale defined thanks to the auxiliary martingales $X_t$ and $Y_t$ as follows
\[
	X_t = Qf(B^M_t,B_t) = Qf(B^M_0,y ) + \int_0^t (\nabla,\partial_y) Qf(B^M_s,B_s) d(U_s dW_s, B_s),
\]

\[
	Y_t =  \int_0^t  \nabla Qf(B^M_s,B_s) dB_s,
\]

\[
	Z_t =  M_t \int_0^t M_s^{-1} dY_s,
\]
where $Y_t$ is by construction differentially subordinate to $X_t$.

\subsection{Main results}

We prove in Theorem \ref{unweightedRiesz} a dimensionless estimate in $L^p$ spaces for the Riesz vector on manifold with non-negative curvature. The first proofs of this result are recent \cite{CD}, \cite{LiErratum}, \cite{BO} and all based on a form of a Bellman function. Our proof is via a sparse domination with continuous index. All these cited Bellman proofs give a better numeric estimate than our proof, but as mentioned earlier, our proof extends (for free) to the weighted case, which the previous ones do not. Our estimate is linear in $p$, which means proportional to $(p-1)^{-1}$ when $p<2$ and to $p-1$ when $p>2$. We note that \cite{BO} have the best numeric constant in this case. We note also that the proof in \cite{CD} gives the linear estimate with $p$ also in the case where the curvature is merely bounded below (and possibly negative) with an appropriately defined Riesz vector involving a Laplacian with a modified spectrum. We do not pursue this here, although parts of our arguments clearly go through also in this case.

\begin{theorem}[$L^p$ estimate]\label{unweightedRiesz}
Suppose that $M$ is a complete Riemannian manifold without boundary and $\Ric_{\varphi}\geq0$. Then for all $f \in C_c^{\infty}(M)$, $p\in (1,\infty)$, we have the following dimension-free estimate

\begin{equation} \label{E1}
	\| R f \|_{L^p(T^*_xM)} \leq 32 \dfrac{p^2}{p-1}\|f\|_{L^p(M)}.
\end{equation}

\end{theorem}

We prove also a dimensionless weighted estimate in $L^p$ spaces for the Riesz vector on manifold with non-negative curvature. In the Euclidean setting, see \cite{DPW}. For the case of manifolds, such an estimate was previously only known in the case $p=2$ see \cite{D}. A priori the weight has to be globally in $L^2$ so as to be able to define the flow characteristic. 
$$\tilde{Q}_p(w)=\sup_{x,y}(Q(w))(x,y)(Q(w^{-\frac1{p-1}}))^{p-1}(x,y).$$
The collection of weights for which this characteristic is finite is denoted $\tilde{A}_p$.
There is also a natural way to extend the class of the weights to resemble more the classical case allowing local $L^1$ weights. In this case we require that constants are integrable in $M$ with the measure $d\mu_{\varphi}$ so as to prove the theorem for cut weights, such as in \cite{D}, that are in $L^1 \cap L^{\infty} \cap L^2$ and then define the characteristic by a limiting procedure and deduce the theorem. See \cite{D} for detailed exposition in the case $p=2$.

\begin{theorem}[weighted $L^p$ estimate]\label{Riesz}
Suppose that $M$ is a complete Riemannian manifold without boundary and $\Ric_{\varphi}\geq0$. Then for all $f \in C_c^{\infty}(M)$, $p\in (1,\infty)$ and $ w \in \widetilde{A}_p$, we have the following dimension-free estimate

\begin{equation} \label{E1w}
	\| R f \|_{L^p(T^*_xM, w)} \leq 32 \dfrac{p^2}{p-1}\widetilde{Q}_p(w)^{\max(1,\frac{1}{p-1})}\|f\|_{L^p(M, w)}.
\end{equation}

\end{theorem}

The technique used in this paper resembles the sparse domination principle for discrete time martingale transforms which originally appeared in \cite{Lacey2015}. This technique has witnessed considerable efforts in the last several years and has been used to prove numerous new results in harmonic analysis, using sparse operators defined on cubes. These cannot give dimensionless estimates, nor are satisfactory results known in the non-doubling case. As in \cite{DP2} we use a sparse operator with continuous stopping times, dominating Li's process $Z_t$ whose projection is the Riesz vector. This is what enables us to use the flow itself without cutting it into cubes, thus resulting in clean dimensionless estimates.

Following \cite{DP2}, we say that the operator $X \mapsto S(X)$ is called sparse if there exists an increasing sequence of adapted stopping times $0=T^{-1}\leq T^0 \leq \cdots $ with nested sets $E_j = \{T^j < \infty\}$, $E_j \subset E_{j-1}$ so that
\begin{equation}\label{sparse1}
S(X) = \sum_{j=-1}^{\infty}X_{T^j} \chi_{E_j} {\text{ where }} X_{T^j} = \mathbb{E}(X|\mathcal{F}_{T^j} );
\end{equation}
\begin{equation}\label{sparse2}
\forall A_j \subset E_j, \ A_j \in \mathcal{F}_{T_j} {\text{   there holds }} \mathbb{P}(A_j \cap E_{j+1}) \leq \dfrac{1}{2}\mathbb{P}(A_j).
\end{equation}

The estimate we aim to show will be a consequence of a sparse domination of the stochastic process $Z_t$ (see [NL,L,DP]). Other than in \cite{DP2} the object is not a martingale, so the sparse domination is different and the key of the proof relies on the weak-$L^1$ estimate for the maximal function of the studied stochastic operator.
We do not aim at the fullest generality here, keeping our goal in mind, an estimate of the Riesz vector. Certain assumptions can certainly be weakened, as the attentive reader will observe.

\begin{lemma}[Weak-type estimate]\label{L: weak type}
Let $X$ be a real valued continuous path martingale and $Y$ a vector valued continuous path martingale so that $Y$ is differentially subordinate with respect to $X$. Let further $Z$ a continuous path semi-martingale whose increments satisfy $dZ_t=V_tZ_tdt +dY_t$ with $V_t$ continuous adapted process with values in non-positive, symmetric $d \times d$ matrices. Let $\lambda >0$. We have
\begin{equation*}
\mathbb{P}\left( (|Z_t|+|X_t|)^*  \geq  \lambda \right) \leq 2 \lambda ^{-1} \|X\|_1.
\end{equation*}
\end{lemma}

\begin{theorem}[Sparse domination]\label{T: sparse decomposition}
Let $X$ be a real valued non-negative continuous path martingale and $Y$ a vector valued continuous path martingale so that $Y$ is differentially subordinate with respect to $X$. Let further $Z$ a continuous path semi-martingale whose increments satisfy $dZ_t=V_tZ_tdt +dY_t$ with $V_t$ continuous adapted process with values in non-positive, symmetric $d \times d$ matrices. Then there exists a sparse domination such that
\[Z^* \leq 8S(X).\]
\end{theorem}

where we recall that we denote by $Z^*=\sup_{t\geq 0}|Z_t|$ the maximal function associated with $Z$.

\begin{theorem}[Weighted estimate]\label{T: weighted estimate for Z}
Let $X$ be a real valued non-negative continuous path martingale and $Y$ a vector valued continuous path martingale so that $Y$ is differentially subordinate with respect to $X$. Let further $Z$ a continuous path semi-martingale whose increments satisfy $dZ_t=V_tZ_tdt +dY_t$ with $V_t$ continuous adapted process with values in non-positive, symmetric $d \times d$ matrices. Then there holds the weighted estimate 
\[ \|Z^{*}\|_{L^p(w)} \lesssim \Phi_p(Q^{\mathcal{F}}_p(w)) \|X\|_{L^p(w)},\] 
where $\Phi_p(x)=x^{\max\{1, \frac1{p-1}\}}$.

\end{theorem}

In general for filtered spaces, the $A_p$ characteristic of $w$ (identified with its closure) is 
$$Q_p^{\mathcal{F}}(w)=\sup_{\tau}\| \mathbb{E} ( (\frac{w_{\tau}}{w})^{\frac1{p-1}} \mid \mathcal{F}_{\tau} )^{p-1}  \|_{\infty}.$$
In the case of interest to us, the characteristic that appears is the one that corresponds to the filtration used by Li at height $y$, denoted $\mathcal{F}^{(y)}$. 
It can be seen, similarly as is known to the Euclidean case, that these characteristic, in a limiting sense, is comparable to the Poisson flow characteristic.

\section{The stochastic process $Z$}
In this section, we prove Lemma \ref{L: weak type} and  Theorem \ref{T: sparse decomposition}.

\begin{proof} (of Lemma \ref{L: weak type}).

This proof is modelled after the exposition in Wang \cite{Wang}. We aim to show
\begin{equation} \label{weak}
\mathbb{P}\left( (|Z_t|+|X_t|)^*  \geq  \lambda \right) \leq 2 \lambda ^{-1} \|X\|_1.
\end{equation}
Indeed, it suffices to show the inequality for $\lambda =1$. To do this, define functions $V,U:\mathbb{R}\times \mathbb{R}^n \to \mathbb{R}$ by 
$$
V(x,y) = \left\{
    \begin{array}{ll}
        -2|x| & \mbox{when } |x|+|y|<1, \\
        1-2|x| & \mbox{when } |x|+|y| \geq 1.
    \end{array}
\right.
$$

$$
U(x,y) = \left\{
    \begin{array}{ll}
        |y|^2-|x|^2 & \mbox{when } |x|+|y|<1, \\
        1-2|x| & \mbox{when } |x|+|y| \geq 1.
    \end{array}
\right.
$$
Let us first observe that everywhere $V\le U$.
Define the stopping time $$T=\inf \{t \geq 0: |X_t|+|Z_t| \geq 1\}.$$ Then $ |X_T|+|Z_T| \geq 1$ and $ |X_t|+|Z_t| < 1$ for $t<T$.\\
We aim to prove that $\mathbb{E}U(X_T,Z_T) \leq 0$, since $V\le U$ the result will follow (see the end of the argument, where we detail the step). We split $$\mathbb{E}U(X_T,Z_T) = \mathbb{E}(U(X_T,Z_T)\chi_{\{T>0\}})+\mathbb{E}(U(X_T,Z_T)\chi_{\{T=0\}})$$ and we show that these contributions are both non-positive.\\
\\
\textbf{Part 1:} $\{T=0\}$.\\
For such $\omega$ where $T=0$ then by definition of $T$ we have $|X_0|+|Z_0| \geq 1$ and $U(X_0,Z_0)=1-2|X_0|$. Assuming that $|Z_0| \leq |X_0|$, then 
\[1 \leq |X_0|+|Z_0| \leq 2 |X_0|,\]
i.e. $1-2|X_0| \leq 0$ and hence
\[\mathbb{E}(U(X_T,Z_T)\chi_{\{T=0\}})=\mathbb{E}(U(X_0,Z_0)\chi_{\{T=0\}}) \leq 0.\]
\\
\textbf{Part 2:} $\{T>0\}$.\\
By simple calculations on the derivatives of $U$ we check that 
\begin{eqnarray}
\partial_{y_i} U(x,y) &=&  2y_i  \label{signU}\\
\partial^2_{xx} U(x,y)&=& -2 \label{signU2},\\
\partial^2_{xy_j} U(x,y)&=& 0,\\
\partial^2_{y_iy_j} U(x,y)&=& 2\delta_{ij} \label{signU3},
\end{eqnarray}
for $|x|+|y| < 1$ and where $\delta_{ij}$ is the Kronecker delta.\\
On $\{T >0\}$, the process evolves in the set $\{(x,y): \ |x|+|y|<1\}$, in the interior of which the function $U$ is twice differentiable, which means that we have the following It\^o formula
\[U(X_T,Z_T)=U(X_0,Z_0)+I_1+\frac{1}{2}I_2,\]
with
\begin{eqnarray*}
I_1 &=& \int_0^T \partial_x U(X_s,Z_s) dX_s + \sum_i \int_0^T \partial_{y_i} U(X_s,Z^{i}_s), dZ^{i}_s \\
I_2 &=& \int_0^T \partial^2_{xx} U(X_s,Z_s), d\langle X,X \rangle _s  +  2\sum_i \int_0^T\partial^2_{xy_i} U(X_s,Z^{i}_s), d\langle X, Z^{i} \rangle _s \\
&& + \sum_i\sum_j \int_0^T \partial^2_{y_iy_j} U(X_s,Z_s), d\langle Z^{i},Z^{j} \rangle _s  .
\end{eqnarray*}
Let's first study $I_1$:\\
Recall that $Z_t$ satisfies the following stochastic differential equation
\begin{equation}\label{dZ}
dZ_t=V_tZ_tdt+dY_t.
\end{equation}
Now if we replace this formula in the expression of $I_1$, we will obtain a local martingale part which is 
\[\int_0^T \partial_x U(X_s,Z_s) dX_s + \int_0^T\langle \partial_y U(X_s,Z_s), dY_s\rangle\]
and a process 
\[ A_T = \int_0^T\langle \partial_y U(X_s,Z_s), V_sZ_s\rangle ds. \]
We may assume that the local martingale is a true martingale without loss of generality and hence its expectation is null. As for the process $A_T$, by (\ref{signU}) we have $$A_T=2\int_0^T \langle Z_s,V_sZ_s\rangle ds\le 0.$$ The non-positivity holds because the integrand is non-positive as well, since $V$ takes values in the class of non-positive matrices. Notice that just like in \cite{D}, the form of the partial derivative of $U$ in the variable $y$ is crucial.\\
\\
Now we deal with $I_2$:\\
By the formulas (\ref{signU2})-(\ref{signU3}), we obtain that
\[\frac12 I_2 = ( \langle Z,Z \rangle_T -|Z_0|^2-  \langle X,X \rangle_T + |X_0|^2 )\chi_{\{T>0\}},\]
and hence it suffices to prove 
\begin{equation} \label{wang}
( \langle Z,Z \rangle_T -|Z_0|^2-  \langle X,X \rangle_T + |X_0|^2 )\chi_{\{T>0\}}\leq 0,
\end{equation}
for any stopping time $T$.  Recall that for all $t$ we have $dZ_t=V_tZ_tdt+dY_t$.
Thus by integrating we have,
\[Z_t - Z_0 =\int_0^tV_sZ_sds + Y_t - Y_0.\]
Taking the quadratic covariance on both sides we obtain 
\begin{eqnarray*} 
\langle Z,Z \rangle_t - |Z_0|^2 &=& \langle Y,Y \rangle_t - |Y_0|^2, \ \ \ \forall t\geq 0 \\
& \leq & \langle X,X \rangle_t - |X_0|^2 \ \text{ by differential subordination}
\end{eqnarray*}
which in turn implies that $\mathbb{E}(I_2) \leq 0$.\\
Finally, $U(X_0,Z_0)=Z_0^2-X^2_0 \leq 0$.  

It remains to show the weak estimate (\ref{weak}):\\
We have $V \le U$ everywhere and $\mathbb{E}U(X_T,Z_T)\le 0$. Therefore 
\begin{eqnarray*}
0 & \geq & \mathbb{E}U(X_T,Z_T)\\
   & \geq & \mathbb{E}V(X_T,Z_T)\\
   & = & \mathbb{E}(V(X_T,Z_T)\chi_{\{|X_T|+|Z_T| \geq 1\}}) + \mathbb{E}(V(X_T,Z_T)\chi_{\{|X_T|+|Z_T| < 1\}}) \\
   & = & \mathbb{E}(-2|X_T|\chi_{\{|X_T|+|Z_T| \geq 1\}}) + \mathbb{E}((1-2|X_T|)\chi_{\{|X_T|+|Z_T| < 1\}}) \\
   & = & \mathbb{P}( |X_T|+|Z_T| \geq 1 ) -2 \mathbb{E}|X_T|,
\end{eqnarray*}
from which we deduce $$\mathbb{P}( (|X_t|+|Z_t|)^{*} \geq 1 )\le 2\| X \|_1$$
and so the lemma is proved. \\ 
\end{proof}

\begin{proof} (of Theorem \ref{T: sparse decomposition}).

Now that we have a weak type result by Lemma \ref{L: weak type}, we are able to use a sparse argument as in \cite{DP2}. Recall for convenience we assumed $X$ non-negative.\\
Consider the processes $Z_t^0=\dfrac{Z_t}{X_0}$ and $Y_t^0=\dfrac{Y_t}{X_0}$ and $X_t^0=\dfrac{X_t}{X_0}$. Applying the result obtained in the first step, we know that the measure of the set
\[E_0=\{ \omega \in \Omega : \max \{ Z^{0*}(\omega), X^{0*}(\omega)\} >4\} \]
is small. Indeed,
\[ |E_0| \leq \dfrac{2}{4} \|X^0\|_1\le \frac12.\]
We can associate $T^{-1}=0$ and a stopping time
\[T^0(\omega)=\inf\{t>0:\max \{ |Z_t^{0}(\omega)|, X_t^{0}(\omega)\} >4\}\]
as the hitting time of the set $L=(4,\infty)$, which is finite in $E_0$, almost surely, by definition.\\
The key of the proof, besides the weak type estimate, relies on recursivity in order to construct a sparse operator.\\
To start, let us suppose we have chosen an increasing stopping time sequence $T^k$. Set for times  $t\in I_k(\omega)=[T^{k-1}(\omega),T^k(\omega)[$ and let us recall that on $I_k$ 
$$Z_t=Z_{T^{k-1}}+\int_{T^{k-1}}^{t}dZ_s$$ so that for all times
\begin{equation}\label{Zsum}
Z_t=\sum_{k=0}^{\infty} Z_t \chi_{t\in I_k}=\sum_{k=0}^{\infty} (Z^{(k)}_t+(Z_t-Z^{(k)}_t)) \chi_{t\in I_k}
\end{equation}
with the $Z^{(k)}_t$ constructed below from $Z_t$ for times in $I_k$ by changing the foot. In order to be more precise, let us first define the martingales $X^{(k)}_t$ and $Y^{(k)}_t$. 
$$ {\text{For }} k>0: \;  X^{(k)}_t=X_{T^{k-1}}+\int_{T^{k-1}}^{\max\{T^{k-1},t\}}dX_s.$$ and 
$$Y^{(0)}_t=Y_0+\int_0^{t} dY_s  {\text{ and if }} k>0: \;  Y^{(k)}_t=\int_{T^{k-1}}^{\max\{T^{k-1},t\}}dY_s.$$
Observe that these are martingales in $\mathcal{F}$ for all $k$ and that $Y^{(k)}$ differentially subordinate 
to $X^{(k)}$. Notice that the processes $$X^k_t=\frac{X^{(k)}_t}{X_{T^{k-1}}} {\text{ and }} Y^k_t=\frac{Y^{(k)}_t}{X_{T^{k-1}}},$$ are also adapted in $\mathcal{F}=(\mathcal{F}_t)_{t\ge 0}$ since at times $t<T^{k-1}$ these processes are constant and hence adapted and at later times the denominator is measurable. Notice that the event $\{T^{k-1}<t\} \in \mathcal{F}_t$ since $T^{k-1}$ is a stopping time.

Now set $$Z^{(0)}_t=Z_0+\int_0^{t}dZ_s$$ and let for $k>0$ $Z^{(k)}_t$ be the process satisfying $Z^{(k)}_t=0$ if $t \le T^{k-1}$ and evolving for $t>T^{k-1}$ according to $$dZ^{(k)}_t=V_tZ^{(k)}_t dt + dY_t$$ with initial condition at time $T^{k-1}$ be set 0. Notice that the so defined process $Z^{(k)}_t$ is adapted to $(\mathcal{F}_t)_{t\ge 0}$ and solves   $dZ^{(k)}_t=V_tZ^{(k)}_t dt + dY^{(k)}_t$ for all times with zero increments for $t<T^{k-1}$. For times $t\ge T^{k-1}$ we know that $W^{(k)}_t=(Z-Z^{(k)})_t$ solves the homogenous equation $$dW^{(k)}_t = V_tW^{(k)}_t dt$$ with initial condition $W^{(k)}_{T^{k-1}}=Z_{T^{k-1}}$. 
Now observe that $$d\langle  W^{(k)}_t,W^{(k)}_t\rangle = 2\langle d W^{(k)}_t, W^{(k)}_t \rangle = \langle V_t W^{(k)}_t,W^{(k)}_t \rangle \le 0$$ because $V_t$ takes values in the non-positive matrices. So we have for $t\ge T^{k-1}$ that $$| Z_t-Z^{(k)}_t |= | W^{(k)}_t | \le |W^{(k)}_{T^{k-1}} | = |Z_{T^{k-1}} |$$
which will give us a control on the error term we induced in the sum (\ref{Zsum}). Using similar arguments as above, we can consider $Z^k_t=\frac{Z^{(k)}_t}{X_{T^{k-1}}}$ and retain these properties, now with respect to martingales $X^k$ and $Y^k$.

We now explain how to choose the sequence of stopping times. Set
\[E_k=\{ \omega \in E_{k-1} : \max \{ Z^{k*}(\omega), X^{k*}(\omega)\} >4\}\]
and its associated stopping time $T^k$ of hitting time. By the above, we know that processes $X^k$, $Y^k$ and $Z^k$ satisfy the assumptions of the weak type estimate and we thus control $|E_k|\le \frac12 \|X^k\chiÑ{E_{k-1}}\|_1\le \frac12|E_{k-1}|$. The second technical assumption of sparse operator in this setting can be shown similarly (see \cite{DP2}). By standard arguments we obtain the pointwise domination

\begin{eqnarray*}
Z^*(\omega) &\leq & 8 \sum_{j=-1}^{\infty}X_{T^j}(\omega) \chi_{E_j}(\omega) \\
&=& 8 S(X)(\omega),
\end{eqnarray*}
by considering $E_{-1}=\Omega $.
\end{proof}

\begin{proof} (of Theorem \ref{T: weighted estimate for Z}).

This follows from the sparse domination and the corresponding estimate for the sparse operator, see \cite{DP2}.
\end{proof}

\section{The Riesz vector}

The proof of the main result now follows standard arguments.

Following Li \cite{Li2008}, recall that
\begin{equation*}
X_t = Qf(X_t,B_t)-Qf(X_0,y).
\end{equation*}
By taking the probabilistic representation of the Riesz transform (\ref{rep}), one can write
\begin{eqnarray*}
\|R  f\|_{L^p(w)}^p &\leq & \lim_{y \rightarrow \infty} 2^p \|Z_{\tau}\|_{L^p(w)}^p \\
& \leq &  \lim_{y \rightarrow \infty} 2^p \|Z^*\|_{L^p(w)}^p \\
& \leq & \lim_{y \rightarrow \infty} (32 \dfrac{p^2}{p-1})^p \Phi_p(Q^{\mathcal{F}^{(y)}}_p(w)) \| X\|_{L^p(w)}^p  \\
& \leq & \lim_{y \rightarrow \infty}(32 \dfrac{p^2}{p-1})^p \Phi_p(Q^{\mathcal{F}^{(y)}}_p(w))\left( \|Qf(B^M_{\tau},B_{\tau})\|_{L^p(w)}^p +  \|Qf(B^M_0,y)\|_{L^p(w)}^p  \right) \\
& \leq & (32 \dfrac{p^2}{p-1})^p  \Phi_p(\tilde{Q}_p(w)) \|f(B^M_{\tau})\|_{L^p(w)}^p \\
& \leq & (32 \dfrac{p^2}{p-1})^p  \Phi_p(\tilde{Q}_p(w)) \|f\|_{L^p(w)}^p ,
\end{eqnarray*}
Notice that sparse domination itself depends upon the used filtration (and hence $y$). Here the norm $\|X\|_{L^p(w)}$ is at $t=\infty$, which is $\tau$ in our stopped processes. We use that $\|Qf(B^M_0,y)\|_{\infty} \to 0$ as $y \to \infty$ and $w \in L^1$ by assumption. $Q^{\mathcal{F}^{(y)}}_p(w)$ is the $A_p$ characteristic that corresponds to the filtration when $B_0=y$ and $ \tilde{Q}_p(w)$ is the Poisson flow characteristic.

\end{document}